\documentclass[11pt]{article}
\usepackage[dvips]{graphics}
\title{ {\bf Robust ${\cal D}$-stability of
uncertain MIMO systems: LMI criteria \footnote{This work is
supported by  National Key Project, National Natural Science
Foundation of China(69925307), and the Research Fund for the
Doctoral Program of Higher Education. Email:
longwang@mech.pku.edu.cn }} }
\date{}

\author{Wang Long \ \ \  Wang Zhizhen \ \ \ Wu Baoyu \\ Department of Mechanics and Engineering Sciences,\\
Center for Systems and Control, Peking University, Beijing 100871,
P.R.China\\ Yu Wensheng\\Institute of Automation, Chinese Academy
of Sciences, Beijing, 100080, P.R.China}
\begin{document}
\maketitle
\thispagestyle{empty}
\begin{minipage}[t]{14.6cm}
{\small {\it Abstract:}   The focal point of this paper is to provide
some simple and efficient criteria to judge the ${\cal D}$-stability of  two families of polynomials, i.e., an
interval multilinear polynomial matrix family and a polytopic
polynomial family. Taking advantage of the uncertain parameter
information, we analyze these two classes of uncertain models
and give some LMI conditions for the robust stability of the two
families.  Two examples illustrate the effectiveness of our results.
}
\end{minipage}
\section{Introduction}
\par\indent

The study of robust stability problems under parameter uncertainties has been pioneered
by the Russian scientist Kharitonov(1978). A rich array of useful results have been
developed over the last twenty years$^{\cite{khar}-\cite{bhat}, \cite{heri}-\cite{mans}}$.
Generally speaking, by dealing directly and effectively with the real parameter
uncertainties in control systems, we can identify apriori the critical subset of the uncertain
parameter set over which stability will be violated. The seminal theorem of
$Kharitonov^{\cite{khar}}$ points out: any real parameter interval polynomial family is
Hurwitz if and only if a special subset (called Kharitonov set) is Hurwitz.
To general uncertain systems, edge theorem gives
a positive answer$^{\cite{bart}}$, which is an one-dimensional test.
\par\indent

Consider the  unity feedback system
with an interval plant and a fixed controller in forward path, its
characteristic polynomial is a multilinear function of certain
interval variables $^{\cite{ hern, koka, long}}$. That is to
say, when all but one variables are fixed, all coefficients of the
polynomial are affine linear in the remaining variable. The collection of all such models is called
 multilinear uncertainties structure.
For an MIMO system, if all relationships between any input and any
output belong to corresponding polytopes (the simplest form
is a line), all admissible models form  a polytopic polynomial
matrix family. In the past one or two decades, there is  continually
growing interest in the robust analysis of matrix$^{\cite{mono,
hern, koka, bial,long}}$. Unlike polynomial case, the vertex result does not hold for
interval matrix. In fact, there does not even exist a result
similar to Edge Theorem for general matrix. Many problems still remain open until now.
Work reported to date shows that a reduced dimensional test
holds$^{\cite{long}}$ for polytopic polynomial matrices. Limited by the complexity of such
problems, several methods have been proposed, such as eigenvalues
estimation, $Lyapunov$ approach, algebraic approach and spectrum
theory, etc.$^{\cite{mono}-\cite{hern}}$. Usually, the algebraic
approach based on Kharitonov theorem is hardly effective and
convenient when used for the robust analysis of matrix directly. $Lyapunov$ approach is an
appealing method developed in the context of $Lyapunov$ theory, which presents in the form of
Linear Matrix Inequalities(LMI). However, until now there is few useful result on robust
test for matrix families under parameter uncertainties.
\par\indent

The purpose of this paper is to address the ${\cal D}$-stability of
interval multilinear polynomial
matrix family and polytopic polynomial family. Even though they are nonlinear problems, we
still can establish some efficient robust stability tests, which
are usually negative for a general  nonlinear case.  Recent work addressed the
${\cal D}$-stability problem for polytope of matrices using
$Lyapunov$ approach$^{\cite{hern}}$  and spectrum theory$^{\cite{mono}}$. In this paper,
we adopt both of two methods to analyze
two classes of uncertain MIMO models and
give several LMI criteria on  robust stability. Our technique not only captures the uncertain parameter
information, but also is easy to use. In the end two examples
demonstrate the effectiveness of our results.
\section{Definitions and Notations}
\par\indent

In this paper, we use the following standard notations and definitions.

\noindent
{\bf Definition 1}\ \
Given an open convex region $D$ in the complex plane, a scalar matrix is
termed $D$-stable, if all its eigenvalues lie in $D$; and a polynomial matrix is termed
$D$-stable, if all roots of its determinant lie in $D$; a matrix
family is termed $D$-stable, if all its members are $D$-stable.\\

\noindent
{\bf Definition 2}\ \
Let $M$ be an arbitrary set, we define $convM$  as the convex hull of
$M$, i.e., the smallest convex set which contains $M$.\\

\noindent
{\bf Definition 3}\ \
For a matrix $A$, its right null-space ${\cal N}_A$ is defined as the space
whose every element $N_A$ satisfies $AN_A=0$. For simplicity,
we denote ${\cal N}_A$  a basis for the right
null-space of $A$.\\

\noindent
{\bf Definition 4}\ \
A polynomial matrix is a matrix with all of its entries being
polynomials; an interval multilinear polynomial matrix is a polynomial
matrix with all of its entries being multilinear dependent on
some interval coefficients; a family of such matrices is called interval
multilinear polynomial matrix family, such as the model in (5); a polytopic polynomial
matrix is a matrix with all of its entries being polytopic polynomials; a family of such matrices
is called polytopic polynomial matrix family, such as the model in (6).\\

\noindent
{\bf Definition 5}\ \
Let ${\cal D}\subset {\cal C}$ be an open convex set of the form
$$
{\cal D}=\left\{s\in{\cal C}:\left[
\begin{array}{c}
1\\
s
\end{array}\right]^\ast B\left[
\begin{array}{c}
1\\
s
\end{array}\right]<0\right\},
$$
where $B$ is an $2\times 2$ matrix and $B^\ast=B$. ${\cal D}$ is  called an
LMI region(\cite{gero}\cite{peau}).\\

\noindent
{\bf Definition 6}\ \
 For every $i\in\{1,\dots, N\}$, $A_i(s)$ is an $n\times n$ polynomial
matrix of the form
\begin{equation}
A_i(s)=A_0^i+A_1^i s+\dots+A_l^i s^l,
\end{equation}
and the $n\times{nl}$ scalar matrix ${\cal A}_i \stackrel{\triangle}{=}(A_0^i,\dots, A_l^i)$
is the coefficient matrix of $A_i(s)$.\\

\noindent
{\bf Definition 7}\ \
For every $i\in\{1,\dots,n\},\; j\in \{1,\dots,n\}$, ${\cal P}_{ij}(s)$
is a polytope of polynomials , i.e.,
\begin{equation}
{\cal P}_{ij}(s)=\left \{\sum_{k=1}^m \lambda_k p_{ij}^{(k)}(s):
\ \  \lambda_k \geq 0 ,\; \sum_{k=1}^m \lambda_k=1
\right \},
\end{equation}
where $p_{ij}^{(1)}(s),\dots, p_{ij}^{(m)}(s) \mbox{ are fixed $k$th-order polynomials}$.
Apparently, the vertex set of every
${\cal P}_{ij}(s)$ is
\begin{equation}
K_{ij}(s)=\{p_{ij}^{(k)}(s),\ \  k=1, \dots, m \}
\end{equation}
and the exposed edge set of every
${\cal P}_{ij}(s)$  is included in the following set
\begin{equation}
E_{ij}(s)=\left\{ \lambda
p_{ij}^{(k)}(s)+(1-\lambda)p_{ij}^{(t)}(s),\; k,t=1, \dots, m,\;
k\not= t,\; \lambda\in [0,1] \right\}.
\end{equation}
{\bf Definition 8}\ \
Denote ${\cal Q}=\left\{(q_1,\dots,q_m)^T:\; q_i\in[q_i^L,
q_i^U]\right\}$.  $a_1(q_1,\dots,q_m),\dots, a_N(q_1,\dots,q_m)$ are
multilinear functions of $q_1,\dots,q_m$.\\

\noindent
{\bf Definition 9}\ \
$$
\begin{array}{c}
R=\left(
\begin{array}{c}
\begin{array}{ccc}
I&&\\
&\ddots&\\
&&I
\end{array}
\begin{array}{c}
O\\
\vdots \\
O
\end{array}\\
\begin{array}{c}
O\\
\vdots \\
O
\end{array}
\begin{array}{ccc}
I&&\\
&\ddots&\\
&&I
\end{array}
\end{array}\right)
\end{array}
$$
 is a $2nl\times(l+1)n$-dimensional scalar matrix and
$I, O$ are $n\times n$ unity, zero
matrices respectively. As usual, $\otimes$ stands for the $Kronecker$ product.

\section{Preliminary results}
\par\indent

 In what follows,  two kinds of uncertain models are considered:
\begin{eqnarray}
{\cal MA}(s)&=& \left\{ \sum_{i=1}^N a_i(q_1, \dots,
q_m)A_i(s),\ \ (q_1,\dots,q_m)^T\in{\cal Q} \right\}\\
{\cal PA}(s)& =& \{(p_{ij}(s))_{n \times n}: p_{ij}(s) \in {\cal P}_{ij}(s),\;
i, j=1, \dots, n \}
\end{eqnarray}
The vertex sets of ${\cal MA}(s)$ and ${\cal PA}(s)$, respectively, are
\begin{equation}
\begin{array}{l}
K_{\cal MA}(s)=\left\{ \sum_{i=1}^N
a_i(q_1,\dots,q_m)A_i(s),\; q_j \in
\{q_j^L, q_j^U\},j=1,\dots,m. \right\}\\
K_{\cal PA}(s)=\{(p_{ij}(s))_{n \times n}: p_{ij}(s) \in K_{ij}(s),\;
i, j=1, \dots, n \}.
\end{array}
\end{equation}
Let $P_n^n$ be the collection
of all permutations of $1, 2, \dots, n$,
and define
\begin{equation}
E_{\cal PA}(s)=\left \{(p_{ij}(s))_{n \times n}:
\begin{array}{l}
p_{k l_k}(s)\in E_{k l_k}(s),\; (l_1, \dots, l_n )\in P_n^n ,\; k=1, \dots, n \\
p_{k i_k}(s) \in K_{k i_k}(s),\; i_k=1, \dots, l_k -1, l_k +1, \dots, n
\end{array}
\right \}
\end{equation}
It is easy to see that, $E_{\cal PA}(s)$ is  a subset of  ${\cal PA}(s)$ produced by taking only
one entry from its exposed edge set in every row/column
and all other entries from their vertex sets.
\par\indent

The  lemma below is due to $Henrion$, et al.\\
 {\bf Lemma 1}\cite{hern}\ \ A polynomial matrix $A_i(s)$ is stable
if there exists a matrix $P_i$ solving the LMI
feasibility problem
$$
{\cal N}_{{\cal A}_i}^\ast R^\ast(B\otimes P_i)R{\cal N}_{{\cal A}_i}<0,\quad P_i=P_i^\ast>0.
$$
Where ${\cal N}_{{\cal A}_i}$ is the right null-space of ${\cal
A}_i$, and $\ast$ denotes the conjugate transpose operator.

\noindent
{\bf Lemma 2}\cite{bart}\ \ (Edge Theorem) Suppose ${\cal D} \subset
{\cal C}$ is a simply-connected region, for any polynomial polytope
$\Omega$ without degree dropping, the root set of $\Omega$
is contained in ${\cal D}$ if and only if the root set of all exposed
edges of $\Omega$ is contained in ${\cal D}$.
\par\indent

Another lemma is on the
${\cal D}$-stability of the  family ${\cal PA}(s)$,\\
{\bf Lemma 3}\cite{long}\ \ ${\cal PA}(s)$ is ${\cal D}$-stable if and only if
$E_{\cal PA}(s)$ is ${\cal D}$-stable.\\
Proof:\ \ For all $A^0(s)\in {\cal PA}(s)$, let $A^0(s)=\left(p_{ij}(s)\right)_{n\times n}$, where
$p_{ij}(s)\in{\cal P}(s)$. For simplicity, we write  $p_{ij}$ for every $p_{ij}(s)\in {\cal P}_{ij}(s)$.
In what follows, we will construct several sets in terms of $A^0(s)$. Let
${\cal A}_k=\left \{A_k(i_1,\dots,i_n;s),\ \ i_1, \dots, i_k \in  \{1, \dots, n \}\right\}(k=0,\dots, n)$, where
$$
A_k(i_1,\dots,i_n;s)=\left(
\begin{array}{rccc}
q_{11}& \dots & q_{1k}& \\
\dots & \dots & \dots & \\
q_{i_1 1}& \dots & q_{i_1 k} & \\
\dots & \dots & \dots &(p_{vt}^0)_{n \times (n-k)} \\
q_{i_k 1}& \dots & q_{i_k k}& \\
\dots & \dots & \dots & \\
q_{n1}& \dots & q_{nk}&
\end{array}
\right)\quad
\begin{array}{c}
q_{lt} \in \left\{
\begin{array}{l}
 E_{i_t t}(s),\ \ l=i_t \\
 K_{i_t t}(s),\ \ l \not=i_t
\end{array}
\right.\\
p_{lv}^0 \mbox{ are entries of } A \\
t=1, \dots, k \quad l=1, \dots, n \\
 v=k+1, \dots, n
\end{array}
$$
It is easy to see that $A^0(s)={\cal A}_0$ and ${\cal A}_k\subset{\cal A}_{k+1}$. In the sequel, we will
prove our statement in two steps:\\
1)\ \ Firstly, we will show that ${\cal A}_n$ is ${\cal D}$-stable if $E_{\cal PA}(s)$ is ${\cal D}$-stable.
By definition, we have
$$
{\cal A}_n=\left\{
\begin{array}{l}
A_n(i_1,\dots,i_n;s)\\
=\left(
\begin{array}{rcl}
q_{11}& \dots & q_{1n}\\
\dots & \dots & \dots \\
q_{n1}& \dots & q_{nn}
\end{array}
\right)
\end{array}
\begin{array}{l}
q_{lt} \in \left\{
\begin{array}{l}
 E_{i_t t}(s), \ \ l=i_t \\
 K_{i_t t}(s), \ \ l \not=i_t
\end{array}
\right.\\
t=1, \dots, n \\
 l=1, \dots, n:\;
\end{array}
i_1, \dots, i_n \in  \{1, \dots, n \}
\right \}.
$$
For all $A_n(i_1,\dots,i_n;s)\in{\cal A}_n$, if $(i_1,\dots,i_n)\in P_n^n$, then $A_n(i_1,\dots,i_n;s)\in E_{\cal PA}(s)$.
Otherwise, there must exist some pair $i_s,i_t$ satisfying $i_s=i_t$. Without loss of generality,
suppose $i_1=i_2=1$, namely
$$A_n(i_1,\dots,i_n;s)=\left(
\begin{array}{rccl}
q_{11}& q_{12}& \dots & q_{1n}\\
\dots & \dots & \dots & \dots \\
q_{n1}& q_{n2}& \dots & q_{nn}
\end{array}
\right)
\qquad q_{11} \in E_{11}(s) \quad q_{12} \in E_{12}(s)
$$
By using Laplace formula on the first row of $A_n(i_1,\dots,i_n;s)$, we have
$$
det(A_n(i_1,\dots,i_n;s))=q_{11}M_{11}+q_{12}M_{12}+\sum_{i=3}^n q_{1i}M_{1i}
$$
where $M_{1i}$ is the algebraic complement of $q_{1i}$. By Lemma 2,
$$
\begin{array}{ccl}
A_n(i_1,\dots,i_n;s) \mbox{ is }{\cal D} \mbox{-stable}& \Leftrightarrow &
q_{11}M_{11}+q_{12}^{0}M_{12}+\sum_{i=3}^n q_{1i}M_{1i} \mbox{ and }\\
& &q_{11}^{0}M_{11}+q_{12}M_{12}+\sum_{i=3}^n
q_{1i}M_{1i}\mbox{ are }{\cal D} \mbox{-stable.}
\end{array}
$$
The corresponding matrices are
$$\left(
\begin{array}{rccl}
q_{11}^{0}& q_{12}& \dots & q_{1n}\\
\dots & \dots & \dots & \dots \\
q_{n1}& q_{n2}& \dots & q_{nn}
\end{array}
\right)
\mbox{  and  }
\left(
\begin{array}{rccl}
q_{11}& q_{12}^{0}& \dots & q_{1n}\\
\dots & \dots & \dots & \dots \\
q_{n1}& q_{n2}& \dots & q_{nn}
\end{array}
\right),
 $$
where $q_{11}^{0} \in K_{11},\;q_{12}^{0} \in K_{12}$. For these two uncertain matrices, if they do not belong to
$E_{PA}(s)$, then there must exist at least two equal indexes. Repeat the same process to them, in the end, we have
$$
E_{\cal PA}(s)\mbox{ is } {\cal D}\mbox{-stable}\Rightarrow {\cal A}_n \mbox{  is  } {\cal D}\mbox{-stable}.
$$
2)\ \ Secondly, for all $A_n(i_1,\dots,i_n;s)\in{\cal A}_n$, by using Laplace formula on the $n$-th column of
$A_n(i_1,\dots,i_n;s)$ and  by Lemma 2, we have
$$
{\cal A}_n \mbox{  is }{\cal D}\mbox{-stable  }\Rightarrow {\cal A}_{n-1}\mbox{  is }{\cal D}\mbox{-stable.}
$$
Continuing this process, we have
$$
{\cal A}_k \mbox{  is }{\cal D}\mbox{-stable  }\Rightarrow {\cal A}_{k-1}\mbox{  is }{\cal D}\mbox{-stable.}
$$
Since $A^0(s)={\cal A}_0$, 1) and 2) imply that $A^0(s)$ is ${\cal D}$-stable. That is to say, if $E_{\cal PA}(s)$ is
${\cal D}$-stable, then, for all $A^0(s)\in {\cal PA}(s)$, $A^0(s)$ is ${\cal D}$-stable. By definition, this means that
${\cal PA}(s)$ is ${\cal D}$-stable. This completes the proof of Sufficiency.
\par\indent

Necessity is obvious because  $E_{\cal PA}(s)$ is a subset of ${\cal PA}(s)$. \ \  $\hfill\diamond$
\par\indent

In this paper, we assume that both ${\cal PA}(s)$ and ${\cal MA}(s)$  have  fixed
degrees.
\section{Main results}
\par\indent

Although overbounding is a bit conservative, it still offers a
powerful tool to solve the stability problem for control systems  with
multilinear uncertainties.  \\
{\bf Theorem 1}\ \ ${\cal MA}(s)\subset conv\{K_{\cal MA}(s)\}$.\\
Proof:\ \ For any $A_0(s)\in{\cal MA}(s)$, by induction, we will show $A_0(s)\in conv\{K_{\cal MA}(s)\}$.
Denote $q=(q_1,\dots, q_t,q_{t+1},\dots, q_m)^T\in {\cal Q}$, where $q_1,\dots,q_t$ are interval parameters
and $q_{t+1},\dots,q_m$ are fixed.
\par\indent

If $t=1$, i.e., $q_1\in[q_1^L,q_1^U]$ and $q_2,\dots,q_m$ are are fixed. Then, in this case,
$$
\begin{array}{l}
{\cal MA}(s)= \left\{ \sum_{i=1}^N a_i(q_1, \dots,
q_m)A_i(s),\ \ q_1\in[q_1^L,q_1^U] \right\}\\
K_{\cal MA}(s)=\left\{ \sum_{i=1}^N
a_i(q_1,\dots,q_m)A_i(s),\; q_1 \in
\{q_1^L, q_1^U\} \right\}.
\end{array}
$$
Since $a_1(q_1,\dots, q_m),\dots,a_N(q_1,\dots, q_m)$ are linear in $q_1$, $\sum_{i=1}^N a_i(q_1,\dots,q_m)A_i(s)$
is also  linear in $q_1$. Clearly,
$A_0(s)\in conv\{K_{\cal MA}(s)\}$.
\par\indent

Assume that the claim holds for $t=k$. When $t=k+1$, we have, in this case,
$$
\begin{array}{l}
 {\cal MA}(s)= \left\{ \sum_{i=1}^N a_i(q_1, \dots,
q_m)A_i(s),
\begin{array}{l}
 q_i\in[q_i^L,q_i^U],i=1,\dots,k+1\\
q_{k+2},\dots,q_m \mbox{  are fixed. }
 \end{array}\right\}\\
K_{\cal MA}(s)=\left\{ \sum_{i=1}^N
a_i(q_1,\dots,q_m)A_i(s),
\begin{array}{l}
 q_i\in\{q_i^L,q_i^U\},i=1,\dots,k+1\\
q_{k+2},\dots,q_m \mbox{  are fixed. }
 \end{array}
\right\}
\end{array}
$$
For all $A_0(s)\in{\cal MA}(s)$, there exists an $m$-dimensional vector
$q^0:=(q_1^0,\dots, q_m^0)^T \in {\cal Q}$, satisfying
$$
A_0(s)=\sum_{i=1}^N a_i(q_1^0,\dots,q_m^0)A_i(s).
$$
Since $a_1(q_1,\dots,q_m),\dots,a_m(q_1,\dots,q_m)$ are linear in $q_{k+1}$, we have
$$
\begin{array}{ll}
a_i(q_1^0,\dots,q_k^0,q_{k+1}^0,\dots,q_m^0)
=&\lambda_0 a_i(q_1^0,\dots,q_k^0,q_{k+1}^L,q_{k+2}^0,\dots,q_m^0)\\
&+(1-\lambda_0)a_i(q_1^0,\dots,q_k^0,q_{k+1}^U,q_{k+2}^0,\dots,q_m^0)
\end{array}
$$
for some $\lambda_0\in[0,1]$.
Therefore,
{\small
$$
\begin{array}{rl}
A_0(s)=&
\sum_{i=1}^N  \left(\lambda_0 a_i(q_1^0,\dots,q_k^0,q_{k+1}^L,q_{k+2}^0,\dots,q_m^0)+
(1-\lambda_0)a_i(q_1^0,\dots,q_k^0,q_{k+1}^U,q_{k+2}^0,\dots,q_m^0)\right)A_i(s)\\
=&\lambda_0\sum_{i=1}^N a_i(q_1^0,\dots,q_k^0,q_{k+1}^L,q_{k+2}^0,\dots,q_m^0)A_i(s)\\
&+(1-\lambda_0)\sum_{i=1}^N a_i(q_1^0,\dots,q_k^0,q_{k+1}^U,q_{k+2}^0,\dots,q_m^0)A_i(s)
\end{array}
$$
}
Since both {$\small\sum_{i=1}^N a_i(q_1^0,\dots,q_k^0,q_{k+1}^L,q_{k+2}^0,\dots,q_m^0)A_i(s)$} and
${\small \sum_{i=1}^N a_i(q_1^0,\dots,q_k^0,q_{k+1}^U,q_{k+2}^0,\dots,q_m^0)A_i(s)}$ belong to the set
$$
conv\left\{ \sum_{i=1}^N
a_i(q_1,\dots,q_m)A_i(s),
\begin{array}{l}
 q_i\in\{q_i^L,q_i^U\},i=1,\dots,k\\
q_{k+1},\dots,q_m \mbox{  are fixed. }
 \end{array}  \right\},
$$
we conclude
$$
A_0(s)\in conv\left\{ \sum_{i=1}^N
a_i(q_1,\dots,q_m)A_i(s),
\begin{array}{l}
 q_i\in\{q_i^L,q_i^U\},i=1,\dots,k+1\\
q_{k+2},\dots,q_m \mbox{  are fixed. }
 \end{array}\right\}.
$$
That is to say, the claim holds also for $t=k+1$. Therefore,  our conclusion is  verified inductively.

\ \  $\hfill\diamond$
\par\indent

For all $B_i(s)\in K_{\cal MA}(s)$, rearrange it as $B_i(s)=B_0^i+B_1^i
s+\dots+B_l^i s^l$. Take ${\cal B}_i=(B_0^i,\dots, B_l^i)$ the coefficient matrix of $B_i(s)$.
It is easy to see
that there exist $2^m$ distinct $B_i(s)$. Hence,
$$
 conv \{K_{\cal A}(s)\}=conv\{B_1(s),\dots,B_{2^m}(s)\}.
$$
With Lemma 1 and Theorem 1, we get an LMI condition for robust
stability of interval multilinear polynomial matrix family
${\cal MA}(s)$: \\
{\bf Theorem 2} \ \ ${\cal MA}(s)$ is robust ${\cal D}$\--stable if
there exist some matrices $P_i=P_i^\ast>0, Q$ solving the LMI feasibility
problem
\begin{equation}
\begin{array}{c}
\left[\begin{array}{c}
R\\
{\cal B}_i
\end{array}
\right]^{\ast}\left[
\begin{array}{cc}
B \otimes P_i& Q\\
Q^{\ast}&0
\end{array}\right]
\left[\begin{array}{c}
R\\
{\cal B}_i
\end{array}\right]
<0,\quad i=1,\dots, 2^m.
\end{array}
\end{equation}
Proof: \ \ For every $A(s)\in {\cal MA}(s)$, by virtue of theorem 1, there exist
$\lambda_1,\dots, \lambda_{2^m}\in [0,1]$ such that $\sum_{i=1}^{2^m}\lambda_i=1$ and
$A(s)=\sum_{i=1}^{2^m}\lambda_i B_i(s)$. Moreover, for all $i\in \{1,\dots, 2^m\}$,
$$
\begin{array}{c}
\left[\begin{array}{c}
R\\
{\cal B}_i
\end{array}
\right]^{\ast}\left[
\begin{array}{cc}
B \otimes P_i& Q\\
Q^{\ast}&0
\end{array}\right]
\left[\begin{array}{c}
R\\
{\cal B}_i
\end{array}\right]
<0
\end{array}
$$
$$
\begin{array}{l}
\Leftrightarrow R^\ast \left(B\otimes P_i\right) R+{\cal B}_i^\ast Q^\ast R+R^\ast Q {\cal B}_i<0\\
\Rightarrow \sum_{i=1}^{2^m}\lambda_i\left(R^\ast (B\otimes P_i) R+{\cal B}_i^\ast Q^\ast R+R^\ast Q {\cal B}_i\right)<0\\
\Leftrightarrow R^\ast \left(B\otimes (\sum_{i=1}^{2^m}\lambda_i P_i)\right) R+
(\sum_{i=1}^{2^m}\lambda_i {\cal B}_i)^\ast Q^\ast R+R^\ast Q
(\sum_{i=1}^{2^m}\lambda_i {\cal B}_i)<0
\end{array}
$$
Multiplying  ${\cal N}_{\cal A}$ from the right and ${\cal N}_{\cal A}^\ast$  from
the left, the inequality becomes
$$
\begin{array}{c}
{\cal N}_{\cal A}^\ast
\left[\begin{array}{c}
R\\
\sum_{i=1}^{2^m}\lambda_i{\cal B}_i
\end{array}
\right]^{\ast}\left[
\begin{array}{cc}
B \otimes (\sum_{i=1}^{2^m}\lambda_i P_i)& Q\\
Q^{\ast}&0
\end{array}\right]
\left[\begin{array}{c}
R\\
\sum_{i=1}^{2^m}\lambda_i{\cal B}_i
\end{array}\right]{\cal N}_{\cal A}
<0
\end{array}
$$
Because of ${\cal N}_{\cal A}={\cal N}_{\sum_{i=1}^{2^m}\lambda_i {\cal B}_i}$, we have
$$
{\cal N}_{\cal A}^\ast R^\ast \left(B \otimes (\sum_{i=1}^{2^m}\lambda_i
P_i)\right) R {\cal N}_{\cal A}<0
$$
From $P_i^\ast=P_i>0$, we have that $\sum_{i=1}^{2^m}\lambda_i
P_i= (\sum_{i=1}^{2^m}\lambda_i P_i)^\ast >0$. Now by the Lemma
1, the conclusion is obvious.  \ \  $\hfill\diamond$

 {\bf Remark 1}\ \ The standpoint of Theorem 2 is to transform stability
problem into a positive real-like condition, and the latter can be
solved using the LMI toolbox.
\par\indent

In regard to Lemma 1 and Lemma 3, we claim that the stability of uncertain family ${\cal PA}(s)$ can be
inferred from whether an LMI condition holds or not for $K_{\cal PA}(s)$. This is shown by the following two theorems.\\
{\bf Theorem 3}\ \ ${\cal PA}(s)$ is ${\cal D}$-stable $\Leftrightarrow conv(K_{\cal PA}(s))$
is ${\cal D}$-stable.\\
Proof: \ \
Sufficiency: We will show $E_{\cal PA}(s)\subset conv(K_{\cal PA}(s))$.
\par\indent

For all $A_1(s)\in E_{\cal PA}(s)$, by the definition of $E_{\cal PA}(s)$,
there exists $(l_1,\dots,l_n)\in P_n^n$ such that
$$
A_1(s)=(p_{ij}(s))_{n \times n}:
\begin{array}{l}
p_{k l_k}(s)\in E_{k l_k}(s),\;k=1, \dots, n,\\
  p_{k i_k}(s) \in K_{k i_k}(s)\; i_k=1, \dots, l_k -1, l_k +1, \dots, n
\end{array}
$$
For $p_{1 l_1}(s)\in E_{1 l_1}(s)$, we know that there exist $m$ real numbers
$\lambda_{11},\dots,\lambda_{1m}\in[0,1]$ satisfying
$$
p_{1 l_1}(s)=\sum_{k=1}^m \lambda_{1k} p_{1 l_1}^{(k)}\mbox{   and  }
\sum_{k=1}^m \lambda_{1k}=1.
$$
Using addition of matrices, we have
$$
A_1(s)=\sum_{t=1}^m \lambda_{1t} A_{1 l_1}^{(t)}
$$
where $A_{1 l_1}^{(t)}$ is the matrix that all its entries coincide
with $A_1(s)$ except one, which lies in the first row and the
$l_1$-th column and equals to $p_{1 l_1}^{(k)}$ for every
$t\in\{1,\dots,m\}$. Thus, for everyone of  $\{A_{1 l_1}^{(1)},\dots,A_{1 l_1}^{(m)}\}$, all of its
entries in the
first row belong to  vertex sets.
\par\indent

For every $A_{1 l_1}^{(t)}\;(t=1,\dots,m)$, applying the same process to $p_{2 l_2}$,
we can find $m$ uncertain matrix families, and for every matrix which
belongs to one of those families, all its entries in the first
and second rows belong to the corresponding vertex sets. Continuing
this procedure, we will get $A_1(s)\in conv(K_{\cal
PA}(s))$.
Therefore, $E_{\cal PA}(s)\subset conv(K_{\cal PA}(s))$.
Then,
$$
\begin{array}{rcl}
conv(K_{\cal PA}(s)) \mbox{ is } {\cal D} \mbox{-stable} &
\Rightarrow& E_{\cal PA}(s)
\mbox{ is } {\cal D} \mbox{-stable}\\
\mbox{ ( by Lemma 3)}&\Leftrightarrow &\mbox{  ${\cal PA}(s)$ is
${\cal D}$-stable}.
\end{array}
$$
\par\indent

Necessity: The relationship between ${\cal PA}(s)$ and $conv(K_{\cal PA}(s))$
can be easily established, thereby Necessity is proved. For all $A_1(s)\in conv(K_{\cal PA}(s))$, there exist $n^2$ numbers
$\lambda_{ij}\in[0,1]$ $\sum_{i,j=1}^n \lambda_{ij}=1$ and $n^2$ matrices
$F_{ij}(s)\in K_{{\cal PA}(s)}$ such that
$$
A_1(s)=\sum_{i,j=1}^n \lambda_{ij} F_{ij}(s).
$$
Denote $F_{ij}(s)=\left(f_{hl}^{ij}(s)\right)_{n\times n}$ with $f_{hl}^{ij}(s) \in K_{hl}(s)$
for all $h,l\in\{1,\dots n\}$. By  addition of matrices,
$$
A_1(s)=\left(\sum_{i,j=1}^n \lambda_{ij} f_{hl}^{ij}(s) \right)_{n\times n}
$$
For every $h\in\{1,\dots ,n\}, l\in\{1,\dots ,n\}$, $\sum_{i,j=1}^n
\lambda_{ij} f_{hl}^{ij}(s)$  still belongs to ${\cal P}_{hl}(s)$
whenever $\lambda_{ij}\in[0,1]$, $\sum_{i,j=1}^n \lambda_{ij}=1$ and $f_{hl}^{ij}(s) \in
K_{hl}(s)$. Thus $conv(K_{\cal PA}(s))\subset {\cal PA}(s)$. This
completes the proof.
\ \ $\hfill\diamond$
\par\indent

For all $A(s)\in {\cal PA}(s)$, we have $A(s)=(p_{ij}(s))_{n\times n}$, where
$p_{ij}(s)\in{\cal P}_{ij}(s)$ with degree  $l$. Rewriting it as
$$
A(s)=A_0+A_1 s+\dots+A_l s^l.
$$
 Denote
${\cal A}\stackrel{\triangle}{=}(A_0,\dots, A_l)$, then ${\cal A}$ is an $n\times {nl}$
scalar matrix.
By Theorem 3, we have\\
{\bf Theorem 4}\ \ ${\cal PA}(s)$ is robust ${\cal D}$-stable if
there exist some matrices $P_{\cal A}=P_{\cal A}^\ast>0, Q$ solving the LMI feasibility
problem
\begin{equation}
\begin{array}{c}
\left[\begin{array}{c}
R\\
{\cal A}
\end{array}
\right]^{\ast}\left[
\begin{array}{cc}
B \otimes P_{\cal A}& Q\\
Q^{\ast}&0
\end{array}\right]
\left[\begin{array}{c}
R\\
{\cal A}
\end{array}\right]
<0,\ \ {\cal A}\in K_{\cal PA}(s).
\end{array}
\end{equation}
Proof:\ \ By  Theorem 3, this problem is equivalent to
the ${\cal D}$-stability of $conv(K_{\cal PA}(s))$. For the latter,
for all $A(s)\in conv(K_{\cal PA}(s))$, there exist $n^2$ numbers
$\lambda_{ij}\in[0,1],\; \sum_{i,j=1}^n \lambda_{ij}=1$ and $n^2$ matrices
$F_{ij}(s)\in K_{{\cal PA}(s)}$ such that
$$
A(s)=\sum_{i,j=1}^n \lambda_{ij} F_{ij}(s).
$$
Now by a similar argument as in the proof of Theorem 2, we get the desired result. \ \ $\hfill\diamond$

\section{Illustrative Examples}
\par\indent

In this section, we give examples to illustrate the utility of our main results. Example 1 is considered
in the context of robust stability of interval multilinear
polynomials with respect to left half plane.\\
{\bf Example 1}\ \ (n=1, N=2, l=3, m=3) Let $A_1(s),A_2(s)$ be two given polynomials
\begin{eqnarray}
A_1(s)=s^3+2.64 s^2+1.82 s+0.37\\
A_2(s)=s^3+5.57 s^2+9.04 s+3.85
\end{eqnarray}
And the uncertain model is ${\cal A}(s)=\left\{a_1(q_1, q_2, q_3)A_1(s)+a_2(q_1, a_2, q_3)A_2(s)\right\}$,
where $q_1\in[1,2],\; q_2\in[3,3.8],\; q_3\in[0.5,0.8]$ and
$
a_1(q_1, q_2, q_3)=0.6 q_1+0.1 q_2-q_3+0.1 q_1 q_2$,
$a_2(q_1, q_2, q_3)= -0.6 q_1-0.1 q_2+q_3+1-0.01 q_2 q_3.
$
For $Hurwitz$ stable, the $2\times 2$ Hermite matrix $B$ corresponds to
$\left[\begin{array}{cc}
0&1\\
1&0
\end{array}\right]$.
Applying Theorem 2 to this problem, it suffices to solve 16 linear matrix inequalities.
Using the LMI Toolbox in $Matlab$, then it is easy to check that the corresponding LMI
problem is feasible. Thus, we conclude that the whole polynomials family is robust
$Hurwitz$ stable.

\noindent
{\bf Example 2}\ \ (n=3, N=2, l=3, m=3) Consider the third order uncertain model below
$$
\begin{array}{l}
{\cal A}(s)=\left\{a_1(q_1, q_2, q_3)A_1(s)+a_2(q_1, a_2, q_3)A_2(s)\right\}\\
\small
 A_1(s)=\frac{1}{10}\\
 \left(
\begin{array}{ccc}
\arraycolsep=1mm
\begin{array}{c}
15 s^3+2.5 s^2+12 s-1\\
10s^2-7.4 s+2.8\\
12 s^2-1.1 s-4.7
\end{array}&
\begin{array}{c}
13 s^3+12 s^2+0.22 s-1.1\\
23 s^3+4.7 s^2-1.3 s+3\\
2.8 s^2+0.99 s-0.5
\end{array}&
\begin{array}{c}
2.7s^3-7.3 s^2-7.9 s-3.7\\
11 s^3-16 s^2+15s+7.6\\
25 s^3-2.2 s^2-3.8 s+0.035
\end{array}
\end{array}\right)\\
\small A_2(s)=\frac{1}{10}\\
 \left(
\begin{array}{ccc}
\arraycolsep=.2mm
\begin{array}{c}
20 s^3-13 s^2-18 s-0.96\\
-9.8 s^2+14 s-4.5\\
6 s^2+17 s+7.1
\end{array}&
\begin{array}{c}
-3.7 s^3+13 s^2+0.078 s-6.7\\
19 s^3+0.35 s^2-7.5 s+3.4\\
-1.3 s^2+1.1 s+1.5
\end{array}&
\begin{array}{c}
-20 s^3-4.5 s^2+9.4 s-0.27\\
-11 s^3-5.5 s^2-20 s+3.8\\
10s^3-1.9 s^2- 10s-5.1
\end{array}
\end{array}\right)
\end{array}
$$
with $q_1\in[1,1.2],\; q_2\in[2.1,2.4],\; q_3\in[1.5,1.8]$ and
$a_1(q_1, q_2, q_3)= q_1- q_2+q_3+0.1 q_1 q_2$, $a_2(q_1, q_2, q_3)= -q_1+ q_2-q_3+1-0.01 q_2 q_3.$
In this Example, the quadratic stability region is the unity circle,
thus the associated matrix is
$\left[\begin{array}{cc}
-1&0\\
0&1
\end{array}\right]$. Solving the corresponding  LMI in Theorem 2,
it is easy to show that the uncertain model is robust $Schur$ stable.\\
{\bf Remark 2}\ \ Our results  can also be verified by the plots
of root loci of the whole polynomials family in following figures.
From the plots of root loci, we can see that our LMI criteria are
not very conservative, and can provide correct, effective
information on robust stability of uncertain systems.
\protect

\section{Conclusion}
\par\indent

In this paper, we have dealt with the performance robustness  of
interval multilinear polynomial matrix families and polytopic polynomial matrix families.
Some computationally tractable and nonconservative sufficient conditions for these two classes of
system models have been obtained.

\end{document}